%% file: tverberg.tex
\documentclass[11pt]{article}
\usepackage{anysize,amssymb,amsmath,amsfonts,amsthm}
\usepackage{graphicx,color}
\usepackage{psfrag,url,paralist}

\title{The Topological Tverberg Theorem\\and Winding Numbers}
\author{\textsc{Torsten Sch\"oneborn}\thanks{This is  a condensed
version of the first author's Diplomarbeit \cite{dipl-Schoeneborn}
at TU Berlin, Inst.\ Math., February 2004, \texttt{arXiv:math.CO/0405393}}
\quad and\quad\setcounter{footnote}{6}%
\textsc{G\"unter M.~Ziegler}\thanks{Partially supported by Deutsche
  Forschungs-Gemeinschaft (DFG), via the \emph{Matheon}
  Research Center ``Mathematics in the Key Technologies'' (FZT86),
  the Research Group ``Algorithms, Structure, Randomness'' (Project ZI 475/3),
  and a Leibniz grant}\\
\normalsize Inst.\ Mathematics, MA 6-2, TU Berlin, D-10623 Berlin, Germany\\
\normalsize \url{{schoeneborn,ziegler}@math.tu-berlin.de}}
\date{Final revised version for JCTA, January 22, 2005}
\newtheorem{theorem}{Theorem}[section]
\newtheorem{corollary}[theorem]{Corollary}
\newtheorem{lemma}[theorem]{Lemma}
\newtheorem{proposition}[theorem]{Proposition}
\newtheorem*{False statement}{False statement}
\newtheorem{conjecture}[theorem]{Conjecture}

\theoremstyle{definition}
\newtheorem{definition}[theorem]{Definition}
\newtheorem{example}[theorem]{Example}

\theoremstyle{remark}

\newtheorem{remark*}{Remark}
 \newcommand{\eps}{\varepsilon}
  \newcommand{\R}{\mathbb{R}}
  \newcommand{\Z}{\mathbb{Z}}
  \newcommand{\df}{\Delta_{(d+1)(q-1)}}
\DeclareMathOperator{\codim}{codim}
\DeclareMathOperator{\dist}{dist}
  \newcommand{\wa}{W_{\neq 0}}
\DeclareMathOperator{\susp}{susp}
\DeclareMathOperator{\ind}{ind}
\begin{document}
\maketitle
\begin{abstract}
  The Topological Tverberg Theorem claims that any continuous map
  of a $(q-1)(d+1)$-simplex to~$\R^d$ identifies points from $q$
  disjoint faces.  (This has been proved for affine maps, for $d\le1$,
  and if $q$ is a prime power, but not yet in general.)
  
  The Topological Tverberg Theorem can be restricted to maps of the
  $d$-skeleton of the simplex.  We further show that it is equivalent
  to a ``Winding Number Conjecture'' that concerns only maps of the
  $(d-1)$-skeleton of a $(q-1)(d+1)$-simplex to~$\R^d$.  ``Many
  Tverberg partitions'' arise if and only if there are ``many
  $q$-winding partitions.''
  
  The $d=2$ case of the Winding Number Conjecture is a problem about
  drawings of the complete graphs $K_{3q-2}$ in the plane.  We
  investigate graphs that are minimal with respect to the winding
  number condition.
\end{abstract}

\section{Introduction}

Our starting point is the following theorem from affine geometry.

\begin{theorem}[Tverberg Theorem]\label{TverbergOriginal}
  Let $d$ and $q$ be positive integers.  Any $(d+1)(q-1)+1$ points in
  $\R^d$ can be partitioned into $q$ disjoint sets whose
  convex hulls have a point in common.
\end{theorem}

This result is from 1966, due to Helge Tverberg \cite{tve66}. 
Today, a number of different proofs are known, including
another one by Tverberg \cite{tve81}. We refer to 
Matou\v{s}ek \cite[Sect.~6.5]{MatousekBZ:BU} for background, for a 
state-of-the-art discussion, and for further references.

By $\Delta_N$ we denote the $N$-dimensional simplex, by 
$\Delta_N^{\le k}$ its $k$-skeleton. Usually we will not distinguish 
between a simplicial complex and its realization. 
One can express the Tverberg Theorem in terms of a linear map:

\begin{theorem}[Tverberg Theorem; equivalent version I]%
\label{TverbergEquivalent1} 
For every linear map 
\[
f:\Delta_{(d+1)(q-1)}\rightarrow \R^d
\]
there are $q$ disjoint faces of $\df$ such that their images have a
point in common.
\end{theorem}

\begin{definition}[Tverberg partitions; Tverberg points]
For $d\ge0$ and $k\ge0$, let
\[
f:\df^{\le k}\ \longrightarrow\ \R^d
\]
be a map.
A set $S$ of $q$ disjoint faces $\sigma$ of $\df^{\le k}$ is a 
\emph{Tverberg partition} 
for the map $f$ if the images of the faces in $S$ have a point 
in common, that is, if 
\[
\bigcap_{\sigma\in S}f(\sigma)\neq\emptyset.
\]
Every point in this nonempty intersection is called a 
\emph{Tverberg point} for the map~$f$.
\end{definition}

In terms of this definition, Tverberg's theorem 
has the following brief statement:

\begin{theorem}[Tverberg Theorem; equivalent version II]
For every linear map 
\[
f:\Delta_{(d+1)(q-1)}\ \longrightarrow\  \R^d
\]
there is a Tverberg partition.
\end{theorem}

The ``Topological Tverberg Theorem'' refers to the
validity of this statement for the greater
generality of \emph{continuous} maps~$f$.
Furthermore, ``Sierksma's dutch cheese problem''
asks for the minimal number of Tverberg partions,
for given $d\ge1$ and $q\ge2$.

\subsection{The Topological Tverberg Theorem}

\begin{conjecture}[``Topological Tverberg Theorem'']\label{ttt}
For every continuous map 
\[
f:\Delta_{(d+1)(q-1)}\ \longrightarrow\  \R^d
\]
there is a Tverberg partition.
\end{conjecture}

For $d=0$ this is trivial. For $d=1$ it
is equivalent to the mean value theorem for 
continuous functions $f:\R\to\R$: The median point is a Tverberg point.

For prime $q$ (and arbitrary $d$) the conjecture was first established by
B\'ar\'any, Shlosman and Sz\H{u}cs \cite{bss81}, using 
deleted products. A proof using deleted joins and the 
$\mathbb{Z}_q$-index is given in \cite{MatousekBZ:BU}.
For prime powers $q$ the conjecture was 
first proved by \"Ozaydin \cite{Ozaydin}; different proofs are 
Volovikov \cite{vol96} and 
Sarkaria \cite{Sarkaria-primepower} (see de Longueville \cite{deL01}).
Thus the above conjecture, which has been \emph{proved} only for
prime powers $q$, is known as the ``Topological Tverberg Theorem.''

Furthermore, it is known (and will be used below) that lower 
dimensional cases follow from higher dimensional ones:

\begin{proposition}[de Longueville {\cite[Prop.~2.5]{deL01}}]%
\label{tttFromDToD-1}%
  If the Topological Tverberg Theorem holds for $q$ and $d$, then it
  also holds for $q$ and $d-1$.
\end{proposition}

All cases with $d\ge2$ and non-primepower $q$ remain open. Thus
the smallest unresolved case is  $d=2$, $q=6$: It deals with maps of the
15-dimensional simplex to $\R^2$. 
Below, we will reduce this case to a question about the planar drawings
of the complete graph $K_{16}$.

According to Matou{\v{s}}ek,
``the  validity of the Topological Tverberg Theorem 
for arbitrary (nonprime) $q$ is one of the most challenging problems 
in this field'' (Topological Combinatorics) \cite[p.154]{MatousekBZ:BU}. 

\subsection{Reduction to the {$(d-1)$}-skeleton} 

The classical version of the Topological Tverberg Theorem
deals with maps from the 
entire $(d+1)(q-1)$-dimensional simplex to $\R^d$.
It may seem quite obvious that this can be reduced to the
$d$-skeleton: In a case of ``general position'' 
Tverberg partitions can only involve faces of dimension at most~$d$.
We establish this in Proposition~\ref{ttttodsc}.

The main result of our paper is a reduction one step further:
We prove that the Topological Tverberg Theorem
is equivalent to the ``Winding Number Conjecture,'' Conjecture \ref{wnc},
which concerns maps of the $(d-1)$-skeleton of $\df$.

However, the ``obvious idea'' for a proof does not quite work;
consequently, the equivalence is not necessarily valid
on a dimension-by-dimension basis.
In particular, although the $d=2$ case of the 
Topological Tverberg Theorem would clearly imply the
validity of the Winding Number Conjecture for $d=2$,
we do not prove the converse implication.

\begin{definition}[Winding number of \boldmath$f$ with respect to a point]
Let $f:S^{d-1}\rightarrow\R^d$ be a continuous map and let $p$ be a 
point in $\R^d$. If $f$ does not attain $p$, 
then $f$ defines a (singular) cycle $[f]$ in the reduced homology 
group $\widetilde H_{d-1}(\R^d\setminus\{p\};\Z)\cong\Z$, and
thus we define the \emph{winding number of~$f$} with respect to~$p$ as
\[
W(f,p)\ :=\ [f]\in\mathbb{Z}.
\]
The sign of $W(f,p)$ depends on the orientation of~$S^{d-1}$
and of~$\R^{d}$, but the expression 
\begin{center}
  ``$W(f,x)\ =\ 0$''
\end{center}
is independent of this choice. 
In particular, for $d=1$ we get that $W(f,p)$ is zero if the 
two points $f(S^0)$ lie in the same component of 
$\R\setminus \{p\}$. Otherwise we say that $W(f,p)\neq 0$.

For any $d$-simplex $\Delta_d$ we have 
$\partial\Delta_d=\Delta_d^{\le d-1}\cong S^{d-1}$;
the winding number $W(f,x)$ for maps 
$f:\partial\Delta_d\rightarrow\R^d$ 
and points $x\notin f(\partial\Delta_d)$ is defined the same way.
Again it is well-defined up to~sign, so the condition
``$W(f,x)=0$'' is independent of orientations.
\end{definition}

\begin{conjecture}[Winding Number Conjecture]\label{wnc}
For any positive integers $d$ and $q$ and every continuous map 
$f:\Delta_{(d+1)(q-1)}^{\le d-1}\rightarrow \R^d$ there are $q$ 
disjoint faces $\sigma_1,\dots,\sigma_q$ of 
$\Delta_{(d+1)(q-1)}^{\le d}$ and a point $p\in\R^d$ such that for 
each~$i$, one of the following holds:
\begin{compactitem}[~$\bullet$~]
 \item $\dim(\sigma_i)\le d-1$ and $p\in f(\sigma_i)$,
 \item $\dim(\sigma_i)= d$, and either $p\in f(\partial\sigma_i)$, 
 or $p\notin f(\partial\sigma_i)$ and $W(f|_{\partial\sigma_i},p)\neq  0$.
\end{compactitem}
\end{conjecture}

A set $S=\{\sigma_1,\dots,\sigma_q\}$ of faces for which some
$p$ satisfies the conditions of the Winding Number Conjecture 
will be referred to as a \emph{winding partition}; 
$p$ will be called a \emph{winding point}.

\begin{figure}[ht]
\begin{center}
\includegraphics*[bb= 110 550 470 670,scale=.57]{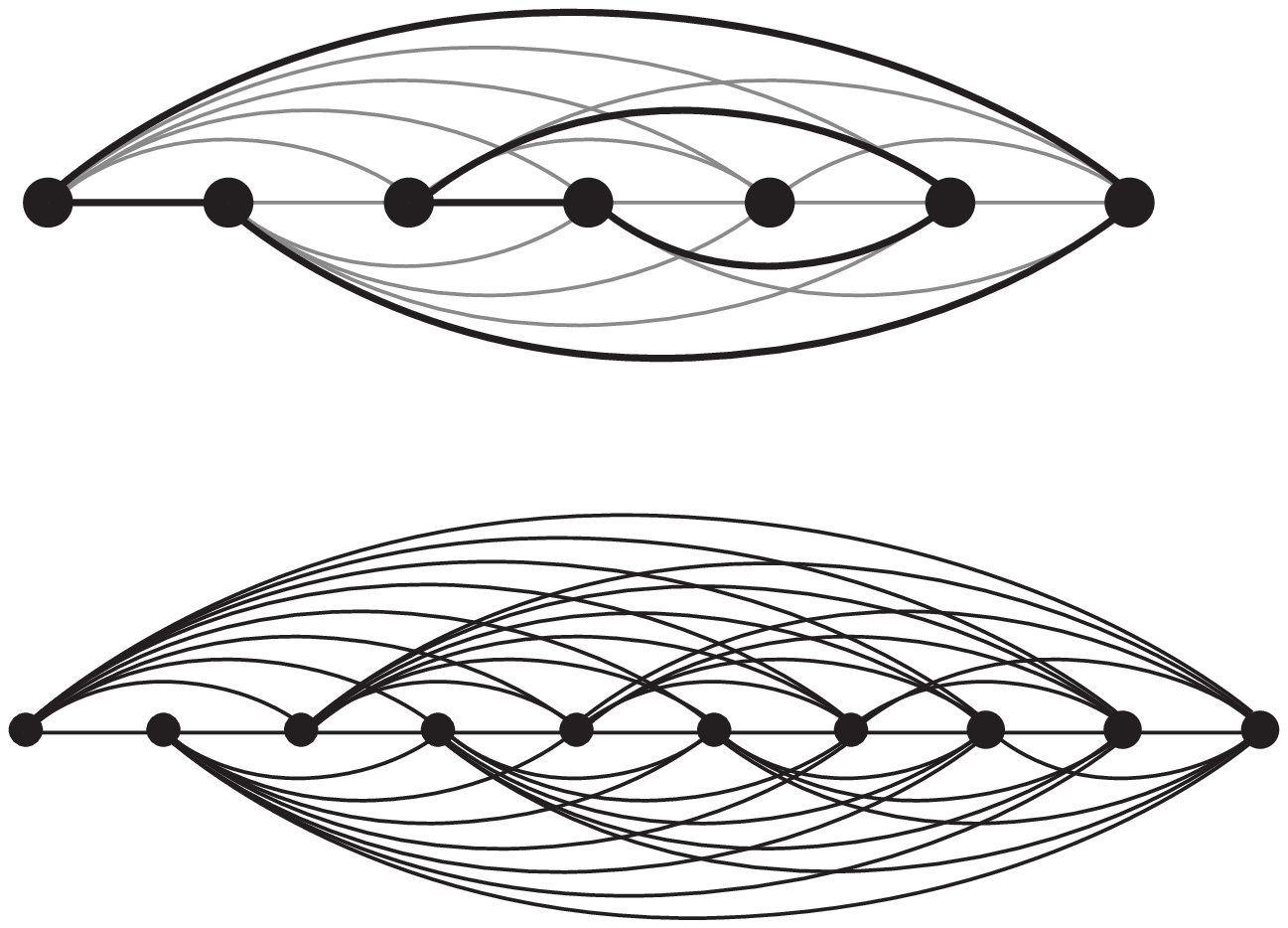}
\includegraphics*[bb= 100 390 500 530,scale=.57]{K_nAlsLinieMitZweiKanten.eps}
\end{center}
\caption{\small The alternating linear drawings of $K_7$ and $K_{10}$. The 
thick edges in the drawing of $K_7$ form a winding partition.} 
\label{K_n}
\end{figure}

\begin{example}\label{K_nByGuy}
  In the case $d=2$, the continuous map $\df^{\le d-1}\to\R^d$ is
  really a drawing of $K_{3(q-1)+1}$, the complete graph with
  $3(q-1)+1=3q-2$ vertices.
  In general, such a drawing may be quite degenerate; it need not be
  injective (an embedding), even locally.
   If the drawing is ``in general position''
  (in a way made precise in the next section), then the Winding Number
  Conjecture says that in the drawing of $K_{3q-2}$ 
  \begin{compactitem}
    \item either $q-1$ triangles 
  (that is, drawings of $K_3$ subgraphs) wind around one vertex,
    \item  or $q-2$ triangles wind around the intersection of two edges,
  \end{compactitem}
  with the triangles, the edges and the vertex being pairwise disjoint
  in~$K_{3q-2}$.

For the ``alternating linear'' drawing
of $K_{3q-2}$ (defined in \cite{saa64};
see Figure~\ref{K_n}) the Winding Number Conjecture is satisfied: 
The $(2q{-}1)$st vertex from the left
is a winding point. For example, the $q-1$ disjoint triangles 
$\langle 1,2,3q-2\rangle$, $\langle 3,4,3q-3\rangle$, \ldots,
$\langle 2q-3,2q-2,3q-q\rangle$ wind around it.
(This is not surprising, since the alternating linear
drawing does have a representation with straight edges,
so the Tverberg Theorem applies, and implies the Winding Number Conjecture
for this example.)
\end{example}

\begin{definition}
For any map $f:\partial\Delta_d\to\R^d$, we define
\[
\wa(f)\ \ :=\ \ 
f(\partial\Delta_d)\ \cup\ 
\{x\in\R^d\setminus f(\partial\Delta_d)\,:\, W(f,x)\neq 0\}.
\]
\end{definition}

\begin{remark*}
It will be advantageous that $\wa(f)$ is a closed set 
containing $f(\partial\Delta_d)$, especially in degenerate cases where 
$\{x\in\R^d\setminus f(\partial\Delta_d)\,:\, W(f,x)\neq 0\}$
might be empty. This is why we add 
$f(\partial\Delta_d)$ to the definition of $\wa(f)$ and include 
``$p\in f(\partial\sigma_i)$'' in our wording of the 
Winding Number Conjecture.
\end{remark*}

\begin{conjecture}[Winding Number Conjecture, equivalent version]
For every continuous map 
$f:\Delta_{(d+1)(q-1)}^{\le d-1}\rightarrow\R^d$ 
there are $q$ disjoint faces $\sigma_1,\dots,\sigma_q$ 
of $\Delta_{(d+1)(q-1)}^{\le d}$ such that 
\[
\bigcap_{\dim(\sigma_i)<d} f(\sigma_i)
\ \cap\bigcap_{\dim(\sigma_i)=d}\wa(f|_{\partial\sigma_i})
\ \ \neq\ \ \emptyset.
\]
\end{conjecture}

This conjecture can be proved easily for $d=1$ (see Proposition 
\ref{wncFuerD=1}). Our main result is the following theorem,
to be proved in the next two sections.

\begin{theorem}\label{WNCequivalenttoTTT}
  For each $q\ge2$, 
  the Winding Number Conjecture is equivalent to the Topological
  Tverberg Theorem.
\end{theorem}

\begin{remark*}%
The basic idea rests on the following two speculations.
\begin{compactitem}[~$\bullet$]
\item Let $F:\df\to\R^d$ be a continuous map. Every winding partition for 
$F|_{\df^{\le d-1}}$ is a Tverberg partition for $F$.
\item Let $f:\df^{\le d-1}\to\R^d$ be a continuous map. Then $f$ 
can be extended to a continuous map $F:\df\to\R^d$ such that 
every Tverberg partition for $F$ is a winding partition for~$f$.
\end{compactitem}
The first statement turns out to be true, but the second one 
needs adjustments, as we will see in the course of the proof. 
\end{remark*}

\subsection{How many Tverberg partitions are there?}

Sierksma conjectured that for every linear map $f:\df\to\R^d$ 
there are at least $((q-1)!)^d$ Tverberg partitions. This number is attained 
for the configuration of $d+1$ tight clusters, with $q-1$ points 
each, placed at the vertices of a simplex, and one point in the 
middle.

For $d=1$ the mean value theorem implies Sierksma's conjecture. 
In almost all other cases, Sierksma's conjecture is still 
unresolved at the time of writing.
Nevertheless, for prime powers~$q$, a lower bound is known
(for the prime case compare 
Matou{\v{s}}ek \cite[Theorem~6.6.1]{MatousekBZ:BU}):

\begin{theorem}[Vu{\v{c}}i{\'c} and {\v{Z}}ivaljevi{\'c} \cite{vuz93},
Hell \cite{hell:_tverb}]\label{ManyTverbergPartitions}
If $q=p^r$ is a prime power, then for every continuous map
$f:\df\to\R^d$ there are at least 
\[
\frac{1}{(q-1)!}\Big(\frac{q}{r+1}\Big)^{\lceil(d+1)(q-1)/2\rceil}
\]
Tverberg partitions.
\end{theorem}

In Section~\ref{sec:HowMany} we discuss how such lower bounds
translate into lower bounds for the number of winding partitions,
with special attention to the case $d=2$, where a direct
translation is not possible. 
We show that Sierksma's conjecture is nevertheless equivalent to
a corresponding lower bound conjecture for the number of
winding partions in any map of the $(d-1)$-skeleton
of a $(d+1)(q-1)$-simplex.

\subsection{Minimal $q$-winding graphs}

The Winding Number Conjecture for $d=2$ is a problem
about the drawings of complete graphs $K_{3q-2}$ in the plane;
it asks whether they are ``$q$-winding'' --- see
Definition \ref{defn:q-winding}.
In Section~\ref{sec:q-winding} we characterize the 
$2$-winding graphs as the non-outerplanar ones, so
$K_4$ is minimal $2$-winding. However, we also show that
$K_7$ is not minimal $3$-winding: Its minimal $3$-winding
subgraph is unique, it is $K_7$ minus a maximal matching.
So the complete graphs $K_{3q-2}$ are not minimal $q$-winding in general.

\section{Reduction to the $d$-skeleton}

The object of this section is to verify that the Topological Tverberg
Theorem guarantees the existence of a Tverberg partition in the
$d$-skeleton of $\df$.

\begin{conjecture}[\boldmath$d$-Skeleton Conjecture]\label{conj:d-SkeletonConj}
  For every continuous map 
\[
  f:\df^{\le d}\ \longrightarrow\ \R^d
\] 
there is a Tverberg partition.
\end{conjecture}

\begin{proposition}\label{ttttodsc}
For each $q\ge2$ and $d\ge1$,
the $d$-Skeleton Conjecture is equivalent to the Topological Tverberg Theorem.
\end{proposition}

It is obvious that the $d$-Skeleton Conjecture implies the Topological
Tverberg Theorem. The converse is harder. 
For this, we verify that any map in question may be approximated
by a piecewise linear map in general position, for which
codimension counts yield that only simplices of dimension
at most $d$ can be involved in a Tverberg partition.
(This needs some care with the definition of 
``general position,'' but is rather straightforward otherwise.)

\subsection{Maps in general position}

For the first lemma, we need the following definition.

\begin{definition}[Linear maps; general position]
Let $\Delta$ be a simplicial complex. A map 
$f:\Delta\rightarrow\R^d$ is \emph{linear} if it is linear on 
every face of $\Delta$. Such a linear map $f$ is 
\emph{in general position} if for every set of disjoint faces 
$\{\sigma_1,\sigma_2,\dots,\sigma_q\}$ of $\Delta$ the inequality
\[
\codim(\bigcap^{q}_{i=1}f(\sigma_i))\ \ \ge\ \ 
\sum^q_{i=1}\codim(f(\sigma_i))
\]
holds, where 
$\codim(\tau):=d-\dim(\tau)$ if $\tau\subset\R^d$. We 
use the convention that $\dim(\emptyset)=-\infty$ and thus  
$\codim(\emptyset)=\infty$. Thus in the case
$\bigcap_{i=1}^q f(\sigma_i)=\emptyset$ the general position
condition holds independently of the right hand side, as then
$\codim(\bigcap_{i=1}^q f(\sigma_i))=\infty$.
\end{definition}

\begin{figure}[ht]
\begin{center}
\includegraphics*[bb= 50 570 285 785,scale=.35]{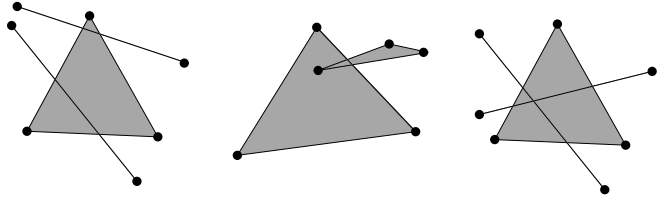}
\includegraphics*[bb= 330 290 570 500, scale=.35]{GeneralPositionExamplesAndCounterexamples.eps}
\includegraphics*[bb= 50 40 300 230, scale=.35]{GeneralPositionExamplesAndCounterexamples.eps}
\end{center}
\vskip-8mm
\caption{\small Images $f(\Delta)$ of linear maps $f:\Delta\to\R^2$ 
in general position.}
\end{figure}

\begin{figure}[ht]
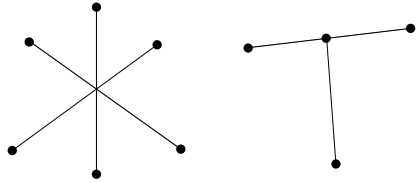

\begin{center}
\includegraphics*[bb= 40 280 300 500, scale=.35]{GeneralPositionExamplesAndCounterexamples.eps}
\includegraphics*[bb= 310 580 560 785,scale=.35]{GeneralPositionExamplesAndCounterexamples.eps}
\end{center}
\vskip-8mm
\caption{\small Images $f(\Delta)$ of linear maps $f:\Delta\to\R^2$ 
\emph{not} in general position. In the last picture, the complex $\Delta$ 
consists of two edges.}
\end{figure}

\begin{definition}[Piecewise linear maps; general position]\label{def2:genpos}
Let $\Delta$ be a simplicial complex. A map 
$f:\Delta\rightarrow\R^d$ is \emph{piecewise linear} if there is a 
subdivision $s:\Delta'\to\Delta$ such that the composition 
$f\circ s:\Delta'\to\R^d$ is a linear map. Furthermore, we
say that $f$ is 
\emph{in general position} if there is a subdivision $s$ such that the 
linear map $f\circ s$ is in general position.
\end{definition}

\begin{figure}[ht]
\begin{center}
\includegraphics*[bb= 40 620 570 730, scale=.65]{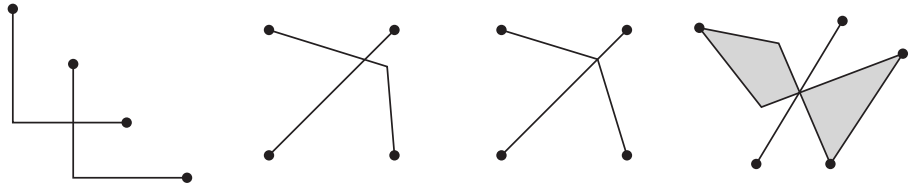}
\end{center}
\vskip-8mm
\caption{\small Images $f(\Delta)$ of piecewise linear maps $f:\Delta\to\R^2$. In 
the first three pictures, $\Delta$ consists of two edges, in the 
last picture $\Delta$ consists of a triangle and an edge. The two 
pictures on the left are in general position, the two on the right are 
not.}\label{PLexamples}
\end{figure}

Whether $f\circ s$ is in general position depends on the subdivision $s$. For 
example, the  map $f$ depicted on the very left in Figure 
\ref{PLexamples} combined with the second barycentric subdivision 
gives a linear map not in general position, although $f$ itself is in general position.

The definition of general position made here may seem overly restrictive for 
the purpose of this section, but we need it in Proposition 
\ref{wnckgeq3}.

\begin{lemma}\label{lemma:codim}
Let $\Delta$ be a simplicial complex and 
$f:\Delta\rightarrow\R^d$ a piecewise linear map in general position. If 
$\{\sigma_1,\sigma_2,\dots,\sigma_q\}$ is a set of disjoint faces 
of $\Delta$, then
\[
\codim(\bigcap^{q}_{i=1} f(\sigma_i))\  \ge\ 
\sum^q_{i=1}\max\{0,d-\dim \sigma_i\}.
\]
\end{lemma}

Here $\bigcap_{i=1}^q f(\sigma_i)$ might have parts of different 
dimension. For polyhedral sets $A$ and $B$ we have
$  \dim(A\cup B)=\max\{  \dim A,  \dim B\}$ and thus
$\codim(A\cup B)=\min\{\codim A,\codim B\}$.

\begin{proof}
Let $s:\Delta'\to\Delta$ be a subdivision such that $f\circ s$ is 
a linear map in general position.
\begin{eqnarray*}
\codim(\bigcap^{q}_{i=1}f(\sigma_i)) & = & 
\min_{\stackrel{\tilde{\sigma}_i\subset\sigma_i}{\tilde{\sigma}_i\mathrm{\ simplex\ in\ }\Delta'}}\codim(\bigcap^{q}_{i=1}f\circ s(\tilde{\sigma}_i))\\
& \ge & \min_{\tilde{\sigma}_i\subset\sigma_i}\sum^q_{i=1}\codim(f\circ s(\tilde{\sigma}_i))\\
& = & \min_{\tilde{\sigma}_i\subset\sigma_i}\sum^q_{i=1}(d-\dim(f\circ s(\tilde{\sigma}_i)))\\
& = & \sum^q_{i=1}(d-\max_{\tilde{\sigma}_i\subset\sigma_i}\dim(f\circ s(\tilde{\sigma}_i)))\\
& \ge & \sum^q_{i=1}(d-\min\{d,\dim\sigma_i\}) 
\ \ =\ \ \sum^q_{i=1}\max\{0,d-\dim\sigma_i\}.
\end{eqnarray*}
\vskip-10mm
\end{proof}

We need an approximation lemma to tackle continuous maps.
For our purposes, a version for finite (compact) simplicial
complexes is sufficient. See \cite[\S16]{Munkres:AT}
for techniques in this context.

\begin{lemma}[Piecewise Linear Approximation Lemma]\label{lemma:PLAL}
Let $\Delta$ be a finite simplicial complex, and
let $f:\Delta\to\R^d$ be a continuous map.
Then for each $\eps>0$ there is a piecewise linear map in general position
$\tilde{f}:\Delta\to\R^d$ with $\|\tilde{f}-f\|_{\infty}<\eps$, with
\[
\|\tilde{f}-f\|_{\infty}\ :=\ \max\{|\tilde{f}(x)-f(x)|:x\in\Delta\}.
\]
\end{lemma}

\subsection{Tverberg partitions in the \emph{d}-skeleton}

\begin{lemma}\label{dscpart1}
  Any Tverberg partition for a piecewise linear map $f:\df\to\R^d$
  in general position contains only faces of dimension at most $d$.
\end{lemma}

\begin{proof}
Let $f$ be in general position, with a Tverberg partition 
$\{\sigma_1,\sigma_2,\dots,\sigma_q\}$. Then
\begin{eqnarray*}
d & \stackrel{(1)}{\ge} & \codim(\bigcap^{q}_{i=1}f(\sigma_i)) \\
& \stackrel{(2)}{\ge} & \sum^q_{i=1}\max(0,(d-\dim\sigma_i))\\
& \stackrel{(3)}{\ge} & \sum^q_{i=1}(d-\dim\sigma_i)\\
& = & qd-(\sum^q_{i=1}((\textrm{number of vertices of 
}\sigma_i)-1))\\
& \ge & qd-(\textrm{(number of vertices of }\df)-q)\\
& = & qd-((d+1)(q-1)+1-q)\ \ =\ \ d.
\end{eqnarray*}
Here (1) holds because $\{\sigma_1,\sigma_2,\dots,\sigma_q\}$ 
is a Tverberg partition and thus $\bigcap^{q}_{i=1}f(\sigma_i)\neq\emptyset$.
(2) holds because $f$ is in general position.
In (3) we have equality only if $d-\dim(\sigma_i)\ge 0$, or 
equivalently, if $\dim(\sigma_i)\le d$ for all $i$, which is what 
we had to prove.
\end{proof}

\begin{lemma}\label{dscpart2}
For every continuous map $f:\df^{\le d}\rightarrow\R^d$ there is an
$\eps_f>0$ such that the following holds: If 
$\tilde{f}:\df^{\le d}\rightarrow\R^d$ is a continuous map with 
$\|\tilde{f}-f\|_{\infty}<\eps_f,$ then every Tverberg partition 
for~$\tilde{f}$ is also a Tverberg partition for~$f$.
\end{lemma}

\begin{proof}
We have to show that for each $S$ that is \emph{not} a Tverberg
partition for~$f$, i.e.\ 
\[
\bigcap_{\sigma\in S}f(\sigma)=\emptyset,
\]
then there is an $\eps_S>0$ such that $S$ is not a Tverberg
partition for any $\tilde{f}$ 
with $\|\tilde{f}-f\|_{\infty}<\eps_S$.
Since there are only finitely many choices for $S$, this
implies the lemma (with $\eps_f:=\min_S\eps_S$).

If $\tilde{f}:\df^{\le d}\rightarrow \R^d$ satisfies 
$\|\tilde{f}-f\|_{\infty}<\eps$, then
\[
\bigcap_{\sigma\in S}\tilde{f}(\sigma)
\ \ \subseteq\ \ 
\bigcap_{\sigma\in S}\big\{x\in\R^d\,:\, 
\dist(x,f(\sigma))\le\eps\big\}\ =:\ C_\eps.
\]
Taking $\eps=\frac1n$ we get a chain
$C_1\supset C_{\frac12}\supset C_{\frac13}\supset\cdots\ $
of compact sets. If all of them are non-empty,
then by compactness also the intersection $\bigcap_n C_{\frac1n}$ is non-empty,
and it would consist of Tverberg points:
$C_0=\bigcap_{\sigma\in S}{f}(\sigma)$.
Thus some $C_{\frac1{n(S)}}$ is empty, and we can take $\eps_S:=\frac1{n(S)}$.
\end{proof}

Thus we have established that 
the Topological Tverberg Theorem implies the $d$-Skeleton Conjecture
(Proposition~\ref{ttttodsc}). Indeed, 
for every map $f:\df^{\le d}\rightarrow\R^d$ there is an $\eps_f$ such that
the following holds:
If a map $\tilde{f}:\df\rightarrow\R^d$ has the property that
its restriction to the $d$-skeleton is $\eps_f$-close to
$f$, then the Tverberg partitions in the $d$-skeleton of~$\tilde{f}$ 
are also Tverberg partitions for $f$ (Lemma~\ref{dscpart2}).
Such a map $\tilde{f}$
may be chosen to be general position piecewise linear (Lemma~\ref{lemma:PLAL}).
So if $\tilde{f}$ has any Tverberg partition,
then this lies in the $d$-skeleton (Lemma~\ref{dscpart1}), and thus yields a
Tverberg partition in the $d$-skeleton for $f$.

\section{Reduction to the $(d-1)$-skeleton}\label{sec:d-1Skeleton}

Now we proceed to prove Theorem~\ref{WNCequivalenttoTTT},
the equivalence of the $d$-Skeleton Conjecture
with the Winding Number Conjecture.

It is quite clear that the
Winding Number Conjecture implies the $d$-Skeleton Conjecture:
Every winding partition is indeed a Tverberg partition.
This rests on the fact that if
$x\in\R^d\setminus f(\partial \Delta_d)$
with $W(f,x)\neq0$, then every extension
of $f$ to $\Delta_d$ must hit~$x$.
(Any map $F:\Delta_d\rightarrow \R^d\setminus\{x\}$ is 
nullhomotopic, so its restriction to $\partial \Delta_d$
has winding number~$0$.)

The proof of the converse is harder. For this we want to show 
that any map 
\[
f:\df^{\le d-1}\ \longrightarrow\ \R^d
\]
can be extended to a map 
\[
F:\df^{\le d}\ \longrightarrow\ \R^d
\]
such that every Tverberg partition for~$F$ is 
a winding partition for~$f$. This would be easy to do if for each 
$d$-dimensional face $\sigma\subseteq\df$, we could arrange that
$F(\sigma)\subseteq\wa(f|_{\partial\sigma})$. However, this is 
not always possible because not every continuous map
$f:S^{d-1}\rightarrow\R^d$ 
is nullhomotopic within $\wa(f)$. For this, we
look at two examples. 

\begin{example}
Let $f:S^1\to\R^2$ be the map
illustrated by Figure \ref{d2counterexample}. The topological 
space $\wa(f)$ is homotopy equivalent to the wedge of two spheres 
$S^1$. The fundamental group $\pi_1(\wa(f))$ is 
$\pi_1(S^1\vee S^1)=\mathbb{Z}\ast\mathbb{Z}$, a
free product. The element 
$[f]\in\pi_1(\wa(f))$ can be written as the nonzero term
$aba^{-1}b^{-1}$ if we choose generators $a,b$ of 
$\mathbb{Z}\ast\mathbb{Z}$ as in the figure. 

If we extend $f$ to $B^2$, then the image covers at 
least one of the two ``holes'' in $\wa(f)$ entirely, which are 
$2$-dimensional sets. There is no one-dimensional subset 
$V\subset\R^2$ such that $f$ is contractible in $\wa(f)\cup V$.

The suspension of this map, $\susp f:S^2\to\R^3$, does
not share this problem: We have 
\[
\wa(\susp f)=\susp \wa(f)\simeq\susp(S^1\vee S^1)=S^2\vee S^2
\] 
again, but this time the homotopy group 
$\pi_2(S^2\vee S^2)$ is not a free product but a direct
sum $\mathbb{Z}\oplus\mathbb{Z}$, so
$[\susp f]=\hat a\hat b\hat a^{-1}\hat b^{-1}=0$ in $\pi_2(\wa(\susp f))$.
\end{example}

\begin{example}
For $d\ge 4$ the homotopy group $\pi_{d-1}(S^{d-2})$ is 
nontrivial; for example, the Hopf map $S^3\to S^2$ is not 
nullhomotopic. Choose such a map 
$f:S^{d-1}\to S^{d-2}$ that is not nullhomotopic. Let 
$i:S^{d-2}\to\R^d$ be an embedding into a
$(d-1)$-dimensional linear subspace of $\R^d$. Then 
$\wa(i\circ f)=i(S^{d-2})$, hence $i\circ f$ can not be 
contracted in $\wa(i\circ f)$.

An important difference between this example and the previous one 
is that here, $i \circ f$ can be contracted within the 
$(d-1)$-dimensional subspace that contains $i(S^{d-2})$;
an extension of the range $\wa(i\circ f)$ to a
$d$-dimensional set is not necessary to make the map 
nullhomotopic.
\end{example}

\begin{figure}
\begin{center}
\input{winding2.pstex_t}
\end{center}
\caption{\small A map $f:S^1\to\R^2$ that is not nullhomotopic 
within $\wa(f)$. The shaded area is $\wa(f)$.} 
\label{d2counterexample}
\end{figure}

Because of the problem illustrated by these examples, 
we take a more technical route.
We need an approximation lemma similar to  Lemma~\ref{dscpart2};
it can be proved along the same lines.

\begin{lemma}\label{wncpart2}
For every continuous map $f:\df^{\le d-1}\rightarrow\R^d$ there is 
an $\eps_f>0$ such that the following holds: If 
$\tilde{f}:\df^{\le d-1}\rightarrow\R^d$ is a continuous map with 
$\|\tilde{f}-f\|_{\infty}<\eps_f$, then every winding partition 
for~$\tilde{f}$ is also a winding partition for~$f$.
\end{lemma}

So, if the Winding Number Conjecture holds for piecewise linear maps
in general position, then it also holds for all continuous maps.

\subsection{The case $d\ge3$}

\begin{definition}[Triangulations of {$\R^d$} in general position]
A \emph{triangulation of $\R^d$} is a simplicial complex 
$\Delta$ with a fixed homeomorphism 
$\|\Delta\|\cong\R^d$
that is linear on each simplex. We do not 
distinguish between a face of the triangulation and the 
corresponding set in $\R^d$.

Triangulations $\Delta_1,\Delta_2,\dots,\Delta_\ell$ of~$\R^d$
are \emph{in general position with respect to each other} if 
\[
\codim\big(\bigcap_{i\in S}\sigma_i\big)\ \ge\ \sum_{i\in S}\codim(\sigma_i)
\]
for every subset $S\subset \{1,\dots,\ell\}$ and faces 
$\sigma_{i}$ of $\Delta_i$.
\end{definition}

\begin{proposition}\label{wnckgeq3}
  For $q\ge2$ and $d\ge3$, the $d$-Skeleton Conjecture
  implies the corresponding case of the Winding Number Conjecture.
\end{proposition}

\begin{proof}
By Lemma~\ref{wncpart2} together with the
Approximation Lemma~\ref{lemma:PLAL}, it suffices to prove the
$d$-dimensional Winding Number Conjecture for the case of
piecewise linear maps $f:\df^{\le d-1}\rightarrow\R^d$:
It suffices to show the existence of winding partitions
for general position piecewise linear maps. 
Our proof consists of three steps:
\begin{compactenum}[~1.\,]
\item 
  For every face $\sigma\subset\df^{\le d}$, 
  choose a triangulation $\Delta_\sigma$ of $\R^d$,
  such that triangulations for disjoint faces are in
  general position with respect to each other. 

\item Extend $f$ to a continuous map 
$F:\df^{\le d}\rightarrow\R^d$
that is ``compatible'' with the $\Delta_\sigma$.

\item By the $d$-Skeleton Conjecture, $F$ has a Tverberg
  partition. Show that every Tverberg partition for~$F$ 
  is a winding partition for~$f$. 
\end{compactenum}
\medskip

\noindent
\emph{Step 1:} 
For each face $\sigma$ of $\df^{\le d}$ choose a
triangulation $\Delta_{\sigma}$ of $\R^d$ such that
the following three conditions are satisfied:
\begin{compactitem}[~--~]
\item For each $\sigma$ of dimension $\dim\sigma\le d-1$,
  the set $f(\sigma)$ is a subset of the
  $\dim(\sigma)$-skeleton of~$\Delta_\sigma$.
  In the case $\dim\sigma=d$, we need that $f(\partial\sigma)$ 
  is a subset of the $(d-1)$-skeleton of~$\Delta_\sigma$.
\item
  If $\sigma_1,\dots,\sigma_\ell$ are disjoint faces of $\df^{\le d}$, then
$\Delta_{\sigma_1},\dots,\Delta_{\sigma_\ell}$ are in general position
with respect to each other. 
(This is possible because $f$ is in general
position.  Here we need the restrictive Definition~\ref{def2:genpos}
of ``general position''!)  
\item
For each $\sigma$ of dimension~$d$, the image $f(\partial\sigma)$
should be contained in a triangulated piecewise linear ball $B_\sigma$ that
is a subcomplex of~$\Delta_\sigma$.
\end{compactitem}
\medskip

\noindent
\emph{Step 2:} Now we extend $f$ to a $d$-face 
$\sigma\subset\df^{\le d}$.
For this let $\tau_1,\dots,\tau_k$ be the 
$d$-faces in $B_\sigma\subseteq\Delta_\sigma$ on which the winding
number of~$f|_{\partial\sigma}$ is zero, that is, the $d$-faces in 
$B_\sigma \cap \overline{(\R^d\setminus\wa(f|_{\partial\sigma}))}$. 
Then we have
\[
\wa(f|_{\partial\sigma})\ \subseteq\ 
B_\sigma\setminus
(\stackrel{\circ}{\tau}_1\cup\dots\cup\stackrel{\circ}{\tau}_k)
\ =:\ 
B_\sigma^\circ.
\]
If we choose a point $x_i$ in each $\tau_i$,
then the set $B_\sigma^\circ$ is a retract of
$B_\sigma\setminus\{x_1,\dots,x_k\}$; so this has the homotopy type
of a wedge of $k$ $(d-1)$-spheres.

\begin{figure}
\begin{center}
\input{d=3proof.pstex_t}
\end{center}
\caption{\small Sketch for the extension problem for
$f|_{\partial\sigma}$ within the triangulation~$\Delta_\sigma$.
Here $B_\sigma^\circ$ would consist of the shaded part together
with the complete $(d-1)$-skeleton of~$B_\sigma$.} 
\label{retractions}
\end{figure}
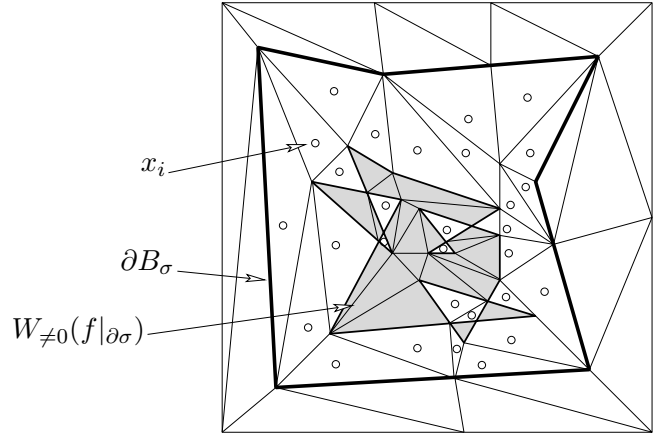

The extension of $f|_{\partial\sigma}$ to $\sigma$ is possible
within $B_\sigma^\circ$
if and only if $f|_{\partial\sigma}$ is 
contractible in $B_\sigma^\circ$, that is,
if the homotopy class
$[f|_{\partial\sigma}]\in\pi_{d-1}(B_\sigma^\circ)$
vanishes. However, for $d\ge3$ we have an isomorphism
\[
\pi_{d-1}(B_\sigma^\circ)\ \cong\ \Z^k.
\]
Furthermore, the boundary spheres of the simplices $\tau_i$
form a homology basis for the wedge, and the evaluation
\begin{eqnarray*}
  \pi_{d-1}(B_\sigma^\circ) &\longrightarrow&\Z^k\\
  f|_{\partial\sigma}      &\longmapsto    &\big(W(f,x_1),\ldots,W(f,x_k)\big)
\end{eqnarray*}
gives this isomorphism.
Thus any map $f|_{\partial\sigma}$ with trivial winding numbers
can be extended to~$\sigma$.
By applying this argument to all $d$-faces of $\df^{\le d}$, 
we obtain a continuous map $F:\df^{\le d}\to\R^d$.
\medskip

\noindent
\emph{Step 3:} We prove that every Tverberg partition for~$F$ is a
winding partition for~$f$. Let $p\in\R^d$ be a Tverberg point and
$\sigma_1,\dots,\sigma_q\subset\df^{\le d}$ a Tverberg partition for~$F$: 
These exist due to the $d$-Skeleton Conjecture for continuous maps.%
\footnote{We have to use the version for continuous maps since
we can not bring a piecewise linear approximation
of $F$ into general position with $F(\sigma)\subset B_\sigma^\circ$.}
 We now show that $\sigma_1,\dots,\sigma_q$ is also a winding
partition for $f$ with winding point $p$:

\begin{compactitem}[~$\bullet$~]
\item $\dim(\sigma_j)\le d-1$: In that case we immediately have 
$p\in F(\sigma_j)=f(\sigma_j)$. 

\item $\dim(\sigma_j)=d$: Suppose $W(f|_{\partial\sigma_j},p)=0$. 
For $1\le i\le q$ let $\tilde{\sigma}_i$ be the face of 
$\Delta_{\sigma_i}$ that contains $p$ in its relative interior, i.e., the 
minimal face containing $p$. We have
\begin{eqnarray*}
d & \stackrel{(1)}{\ge} & \codim(\bigcap^{q}_{i=1}\tilde{\sigma}_i) \\
& \stackrel{(2)}{\ge} & \sum^q_{i=1}\codim(\tilde{\sigma}_i) 
\ \ =\ \ \sum^q_{i=1}(d-\dim(\tilde{\sigma}_i))\\
& \stackrel{(\ast)}{\ge} & \sum^q_{i=1}(d-\dim(\sigma_i))\\
& = & qd-((d+1)(q-1)+1-q)\ \ =\ \ d.
\end{eqnarray*}
where (1) holds because $\bigcap^{k}_{i=1}\tilde{\sigma}_i$ 
contains $p$ and therefore is not empty,
and (2) holds because the $\Delta_{\sigma_i}$ are in general position
with respect to each other.

The inequality $(\ast)$ is an equality if and only if 
$\dim(\tilde{\sigma}_i)=\dim(\sigma_i)$ for all $i$ and in 
particular for $i=j$. Hence 
$\dim(\tilde{\sigma}_j)=\dim(\sigma_j)=d$. Outside of 
$\wa(f|_{\partial\sigma})$, the image $F(\sigma_i)$ lies entirely 
in the $(d-1)$-skeleton of $\Delta_{\sigma_i}$; therefore $p$ 
must lie in $\wa(f|_{\partial\sigma})$.
\end{compactitem}
\vskip-4mm
\end{proof}

\subsection{The case $d=2$}\label{d=2}

We do not know whether the cases $d=2$ of the Winding Number Conjecture
and the $d$-Skeleton Conjecture are equivalent. Thus we take a
different route:

\begin{proposition}\label{wncd2}
For each $q\ge2$, 
if the Winding Number Conjecture holds for $d+1$, then it also holds for $d$.
\end{proposition}

\begin{proof}[Proof (suggested by~{\cite[Prop.~2.5]{deL01};
cf.~Prop.~\ref{tttFromDToD-1}})]
For any continuous map 
\[
f:\Delta_{(d+1)(q-1)}^{\le d-1}\ \longrightarrow\ \R^d
\]
we identify $ \R^d$ with 
$\R^d\times\{0\}=\{x\in\R^{d+1}:x_{d+1}=0\}\subset\R^{d+1}$,
and construct an extension to 
$F:\Delta_{(d+2)(q-1)}^{\le d}\rightarrow\R^{d+1}$, as follows. 
Choose points $Q_1,\dots,Q_{q-1}$ (which may coincide) 
in the upper halfspace $\R^d\times\R^+=\{x\in\R^{d+1}:x_{d+1}>0\}$, 
and an additional cone point $Q$ in the lower halfspace $\R^d\times\R^-$.

We consider $\Delta_{(d+1)(q-1)}^{\le d-1}$ as a subcomplex of 
$\Delta_{(d+2)(q-1)}^{\le d}$, so the latter has $q-1$ additional
vertices $P_1,P_2,\dots,P_{q-1}$.
For $F$, we map the $P_i$ to $Q_i$; 
all faces of $\Delta_{(d+2)(q-1)}^{\le d}$
that involve at least one of the new vertices $P_i$ are 
mapped accordingly by linear extension.
For the $d$-faces of $\Delta_{(d+1)(q-1)}^{\le d}$ we 
perform a stellar subdivision, map the new center vertex 
to~$Q$, and extend canonically.

\begin{figure}[ht]
\begin{center}
\includegraphics[bb = 80 250 530 650, scale=.6,clip]{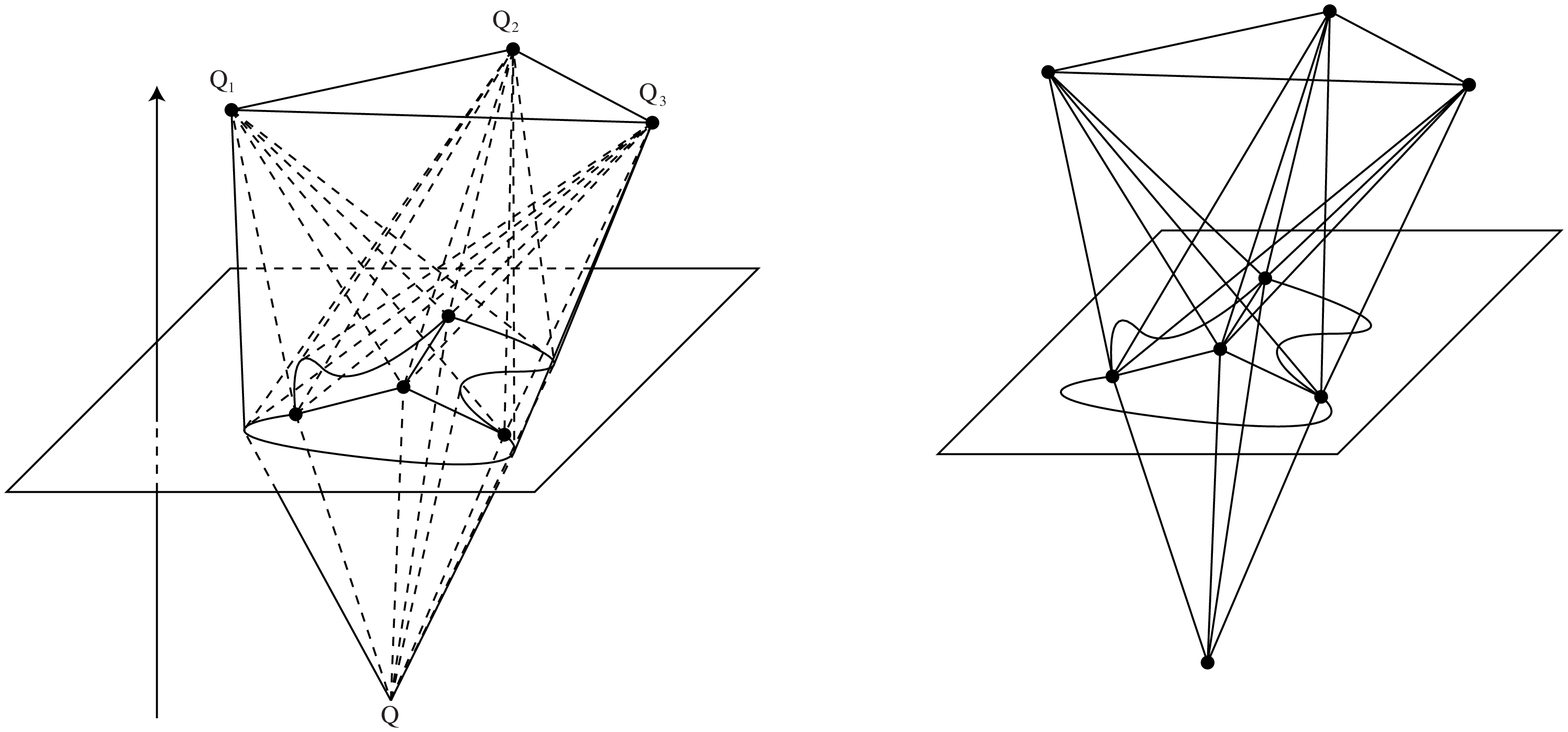} 
\end{center}\vskip-5mm
\caption{\small The map $F$. The plane $\R^2$ contains the image of 
$\Delta_4^{\le1}=K_4$, the three points above the plane are 
$Q_1,Q_2,Q_3$ and the point below is $Q$.} \label{suspension}
\end{figure}

The Winding Number Conjecture for $d+1$ applied to $F$ yields a
winding point $p$ in $\R^{d+1}$ with a winding partition consisting of
$q$ disjoint faces $\sigma_1,\dots,\sigma_q$ of
$\Delta_{(d+2)(q-1)}^{\le d+1}$.

The winding point cannot be in the upper halfspace,
since then all the $F(\sigma_i)$ would need to intersect
the upper halfspace, so the disjoint faces
$\sigma_i$ would need to contain distinct
vertices $P_j$, and there are only $q-1$ of these.
If $p$ were in the lower halfspace, then all the 
$\sigma_i$ would need to be $d$-faces of~$\Delta_{(d+1)(q-1)}^{\le d}$.
For these disjoint $d$-faces we would need
$q(d+1)$ vertices in $\Delta_{(d+1)(q-1)}$, 
which has only $(d+1)(q-1)+1=q(d+1)-d$ vertices.

Thus $p$ has to be in~$\R^d$. 
Define $\tilde{\sigma}_i:=\sigma_i\cap\df^{\le d}$. 
We claim that 
$\tilde{\sigma}_1,\dots,\tilde{\sigma}_q$ are $q$ disjoint faces that
form a winding partition for $f$. For this, we have three cases.
\begin{compactenum}[1.~]
\item If $\dim(\sigma_i)\le d-1$, then
$p\in F(\sigma_i)\cap\R^d=F(\sigma_i\cap\Delta_{(d+1)(q-1)}^{\le
  d-1})=F(\tilde{\sigma}_i)=f(\tilde{\sigma}_i)$.

\item If $\dim(\sigma_i)=d$, then
$p\in F(\sigma_i)\cap\R^d=F(\sigma_i\cap\Delta_{(d+1)(q-1)}^{\le d-1})$.
Now either $\sigma_i\subseteq\Delta_{(d+1)(q-1)}^{\le d}$,
  then we have 
           $F(\sigma_i\cap   \Delta_{(d+1)(q-1)}^{\le d-1})=
             F(\partial\tilde\sigma_i)=f(\partial\tilde\sigma_i)$;
  or   $\sigma_i\not\subseteq\Delta_{(d+1)(q-1)}^{\le d}$,
  so we have 
           $F(\sigma_i\cap   \Delta_{(d+1)(q-1)}^{\le d-1})=
             F(\tilde\sigma_i)=f(\tilde\sigma_i)$.

\item For $\dim(\sigma_i)=d+1$ we may assume that $p$ is 
not in $F(\partial\sigma_i)$. We know that $p$ lies in 
$\wa(F|_{\partial\sigma_i})\cap\R^d$, therefore 
$F(\partial\sigma_i)$ must contain points in both halfspaces. 
Thus $\sigma_i$ contains exactly one of the $P_j$, 
$\tilde{\sigma}_i$ is $d$-dimensional, and
\begin{eqnarray*}
p\ \in\ \wa(F|_{\partial\sigma_i})\cap\R^d & = & \{x\in\R^{d+1}\,:\, 
W(F|_{\partial\sigma_i},x)\neq 
0\}\cap\R^d\\
& =& \{x\in\R^d\,:\, W(f|_{\partial\tilde{\sigma}_i},x)\neq 
0\}\ \ =\ \  \wa(f|_{\partial\tilde{\sigma}_i}).\\[-13mm]
\end{eqnarray*}
\end{compactenum}
\end{proof}

\section{The number of winding partitions and Tverberg partitions}%
\label{sec:HowMany}

The Winding Number Conjecture, and the analogue of the 
Sierksma conjecture for winding partitions, are trivial in the case $d=1$:

\begin{proposition}[The case \emph{d}$=$1]\label{wncFuerD=1}~\\
  For every continuous mapping $f:\Delta_{2(q-1)}^{\le0}\to\R$, there
  are at least $(q-1)!$ winding partitions.
\end{proposition}

\begin{proof}
  $\Delta_{2(q-1)}^{\le0}$ is a set of $2(q-1)+1=2q-1$ vertices.
  $f(\Delta_{2(q-1)}^{\le0})$ is a set of $2(q-1)+1$ real numbers
  (counted with multiplicity). Denote the vertices of
  $\Delta_{2(q-1)}^{\le0}$, ordered by their function value, by
  $P_1,\dots,P_{q-1},M,Q_1,\dots,Q_{q-1}$. A partition of these points
  into $q$ sets is a winding partition for $f$ if one of the sets is
  $\{M\}$ and all the other sets contain exactly one of the $P_i$ and
  one of the $Q_j$. There are $(q-1)!$ such partitions.
\end{proof}

\begin{corollary}[The case \emph{d}$\ge$3]\label{EquivalenceOfLowerBounds}
For each $q\ge2$ and $d\ge3$, the following three numbers are equal:
\begin{compactenum}[\rm1.~]
\item the minimal number of Tverberg partitions for continuous 
maps~$f:\df\to\R^d$,
\item the minimal number of Tverberg partitions for continuous 
maps~$f:\df^{\le d}\to\R^d$,
\item the minimal number of winding  partitions for continuous 
maps~$f:\df^{\le d-1}\to\R^d$.
\end{compactenum}
If $d=2$, then the first two of these numbers are equal.
\end{corollary}

\begin{proof}
By Lemmas \ref{dscpart2} and \ref{wncpart2}, 
the minimal numbers will be achieved 
for general position maps~$f$. For these, all Tverberg partitions
lie in the $d$-skeleton by Lemma \ref{dscpart1}. The
proof of Proposition \ref{wnckgeq3} shows that for $d\ge3$, 
each $f:\df^{\le d-1}\to\R^d$ can be extended to 
an   $F:\df^{\le d}  \to\R^d$ such that all Tverberg partitions
for~$F$ are winding partitions for~$f$.
\end{proof}

If Sierksma's conjecture on the minimal number of Tverberg partitions is 
correct, then the equivalence established in the previous 
proposition carries over to the case $d=2$:

\begin{theorem}
For each $q\ge2$, the following three statements are equivalent:
\begin{compactenum}[\rm1.]
\item \emph{Sierksma's conjecture:} For all positive integers $d$ and 
$q$ and every continuous map $f:\df\to\R^d$ there are at least 
$((q-1)!)^d$ Tverberg partitions.

\item For every continuous map $f{:}\,\df^{\le d}{\to}\R^d$ there are at 
least  $((q-1)!)^d$ Tverberg partitions.

\item For every continuous map $f:\df^{\le d-1}\to\R^d$ there are 
at least $((q-1)!)^d$ winding partitions.
\end{compactenum}
\end{theorem}

\begin{proof}
By our proof for Theorem~\ref{WNCequivalenttoTTT},
we know that Statements 1\ and~2\ are 
equivalent and that Statement 3 implies Statement 2, which in 
turn guarantees Statement~3 if $d\neq 2$.

We now prove that the case $d=3$ of Statement 3 implies the case 
 $d=2$. By Lemma \ref{wncpart2}, it is sufficient to 
examine piecewise linear maps $f:\Delta_{3(q-1)}^{\le1}\to\R^2$ in general position. Regard 
$f$ as a map $\Delta_{3(q-1)}^{\le1}\to\R^3$ in the way we did in 
the proof of Proposition \ref{wncd2}. For each pair $e_1,e_2$ of 
1-dimensional faces of $\Delta_{3(q-1)}^{\le1}$, define one of them to 
be the ``upper'' and the other one to be the ``lower'' one of the 
pair. Now alter $f$ in the following way: For each intersection 
$P\in f(e_1)\cap f(e_2)$ of the images of two edges, change $f$ 
slightly so that the image of the ``upper'' line runs above the 
image of the ``lower'' line at $P$, i.e., has a bigger last 
coordinate (see Figure~\ref{fTofTilde}). We call this new map 
$\tilde{f}:\Delta_{3(q-1)}^{\le1}\to\R^3$.

\begin{figure*}
\begin{center}
\includegraphics[bb = -5 350 640 540, scale = .7,clip]{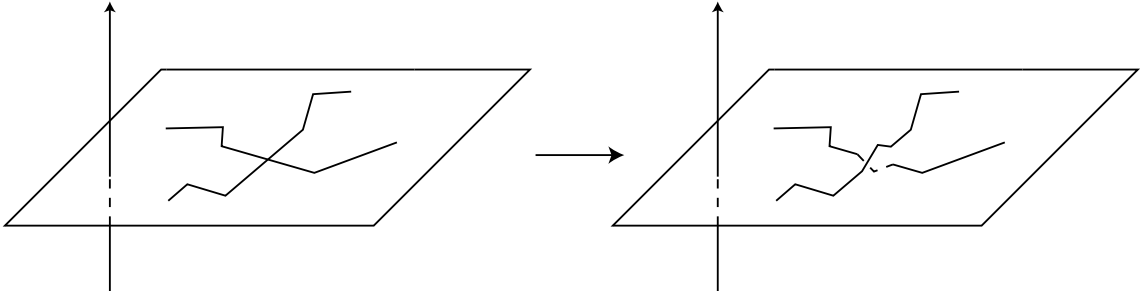} 
\vskip-6mm
\end{center}
\caption{\small How $\tilde{f}$ is obtained from $f$} \label{fTofTilde}
\end{figure*}

We continue similar to the proof of Proposition \ref{wncd2} and 
choose points $Q_1,\dots,Q_{q-1}$ high above $\R^2$ and a 
point $Q$ far below $\R^2$ and extend $\tilde{f}$ to a map 
$F:\Delta_{4(q-1)}^{\le2}\to\R^3$ by taking cones using the $Q_i$ 
and $Q$. Let $\{\sigma_1,\dots,\sigma_q\}$ be a winding partition for $F$ and 
denote $\tilde{\sigma}_i:=\sigma_i\cap\Delta_{3(q-1)}^{\le1}$. By the 
argument given in that proof, 
$\{\tilde{\sigma}_1,\dots,\tilde{\sigma}_q\}$ is a winding partition for $f$. 
Since $f$ is in general position, there are two possibilities for the 
$\tilde{\sigma}_i$.
\begin{compactitem}[~$\bullet$~]
\item 
The $\tilde{\sigma}_i$ are 2-dimensional, except for one,
say $\tilde{\sigma}_1$, that is
0-dimensional. Since $\{\sigma_1,\dots,\sigma_q\}$ is a 
winding partition for our constructed $F$, the faces 
$\sigma_2,\dots,\sigma_q$ have to be 3-dimensional and the face 
$\sigma_1$ has to be 0-dimensional. Therefore each of the faces 
$\sigma_2,\dots\sigma_q$ contains exactly one of the vertices 
$P_i$. Hence the winding partition 
$\{\tilde{\sigma}_1,\dots,\tilde{\sigma}_q\}$ for~$f$ corresponds 
to $(q-1)!$ winding partitions of $F$.
\item 
All but two of the $\tilde{\sigma}_i$ are 2-dimensional,
and the other two, say
$\tilde{\sigma}_1$ and $\tilde{\sigma}_2$, are 
1-dimensional. W.l.o.g.\ let $\tilde{\sigma}_1$ be the ``upper'' 
one of the two. Since $\{\sigma_1,\dots,\sigma_q\}$ is a winding partition 
for~$F$, the faces $\sigma_3,\dots,\sigma_q$ have to be 
3-dimensional, the face $\sigma_2$ has to be 2-dimensional and 
the face $\sigma_1$ has to be 1-dimensional. Hence the winding partition 
$\{\tilde{\sigma}_1,\dots,\tilde{\sigma}_q\}$ for~$f$ corresponds 
to $(q-1)!$ winding partitions of $F$.
\end{compactitem}
In both cases there is a $1$-to-$(q-1)!$ map
between winding partitions of $f$ and
winding partitions of $F$. Since there are at least $((q-1)!)^3$
winding partitions of $F$, there have to be at least $((q-1)!)^2$
winding partitions of $f$.
\end{proof}

\begin{example}
  For the alternating linear model of $K_n$ described in Example
  \ref{K_nByGuy}, there are $((q-1)!)^2$ winding partitions, exactly
  the bound conjectured in the previous Theorem.
\end{example}

We now know that for $d\ge3$ the proved and conjectured lower bounds for the
number of Tverberg partitions also apply to the number of winding
partitions; in the following we derive a nontrivial lower bound 
on the number of winding partitions also for the case $d=2$.

\begin{proposition}[The case \emph{d}$=$2]
Let $q=p^r$ be a prime power. 
For every map $f{:}\,\Delta_{3(q-1)}^{\le1}{\to}\R^2$
there are at least 
\[
\frac{1}{((q-1)!)^2}\Big(\frac{q}{r+1}\Big)^{2(q-1)}
\]
winding partitions.
\end{proposition}

\begin{proof}
In the case $d=3$, there are at least 
$b:=\frac{1}{(q-1)!}\cdot(\frac{q}{r+1})^{2(q-1)}$ 
Tverberg partitions (Theorems 
\ref{ManyTverbergPartitions} and \ref{EquivalenceOfLowerBounds}) 
and thus the same number of winding partitions. By the proof of the 
previous theorem, $\frac{b}{(q-1)!}$ is a bound for the number of 
winding partitions for $d=2$. 
\end{proof}

\section{$q$-Winding  Graphs}\label{sec:q-winding}

For $d=2$ the Winding Number Conjecture claims that complete graphs
$K_{3q-2}=\Delta_{3(q-1)}^{\le1}$ have a certain property. We now
consider all graphs that have this property.
For this, we interpret graphs as ($1$-dimensional)
topological spaces if needed. Thus, a \emph{drawing} of $G$
is just a continuous map $G\rightarrow\R^2$.
(Nevertheless, the \emph{paths} and \emph{cycles} 
in the following definition are required to be subgraphs, 
so the paths start and end at vertices. A single vertex is a path of 
length~$0$.) 

\begin{definition}[{$q$}-winding]\label{defn:q-winding}
A graph is $G$ \emph{$q$-winding} if for every drawing
$f:G\to\R^2$ there are $q$ disjoint paths or cycles 
$P_1,\dots,P_q$ in $G$ with
\[
\Big(\bigcap_{P_i\mathrm{\ is\ a\ path}}
f(P_i)\Big) \cap\Big(\bigcap_{P_i\mathrm{\ is\ a\ cycle}} 
\wa(f|_{P_i})\Big)\ \neq\ \emptyset.
\]
In this situation $P_1,\dots,P_q$ is a \emph{$q$-winding partition} for $f$. 
\end{definition}

The case $d=2$ of the Winding Number Conjecture claims that $K_{3q-2}$
is $q$-winding.
This is proved in the case when $q$ is a prime power.
So the first ``undecided case'' is $q=6$:
Does every drawing of $K_{16}$ have a $6$-winding partition,
into either a vertex and five triangles, or into
two edges and four triangles?

We now take a closer look at 2- and 3-winding  graphs.
(Every non-empty graph is 1-winding.)

\subsection{2-Winding  graphs and $\Delta$-to-$\!Y$ operations}

\begin{proposition}\label{K_4ANDK_2,3}
$K_4$ and $K_{2,3}$ are $2$-winding.
\end{proposition}

Our proof will be phrased in terms of $\Delta$-to-$\!Y$ operations
(compare \cite[Sect.~4.1]{Z35}).
We discuss their effect on $q$-winding graphs in general before 
we return to the proof of the proposition.

\begin{definition}[{\boldmath$\Delta$-to-$\!Y$ operations}]
  A \emph{$\Delta$-to-$\!Y$ operation} transforms a graph $G$ into a 
  graph $G'$ by deletion of
  the three edges of a triangle, and addition of a new
  3-valent vertex that is joined to the three vertices of the
  triangle. A \emph{$Y\!$-to-$\Delta$ operation} is the reverse of a
  $\Delta$-to-$\!Y$ operation.
\end{definition}

\begin{figure}[htb]
\begin{center}
\includegraphics*[bb= 4 560 680 830, scale=.35]{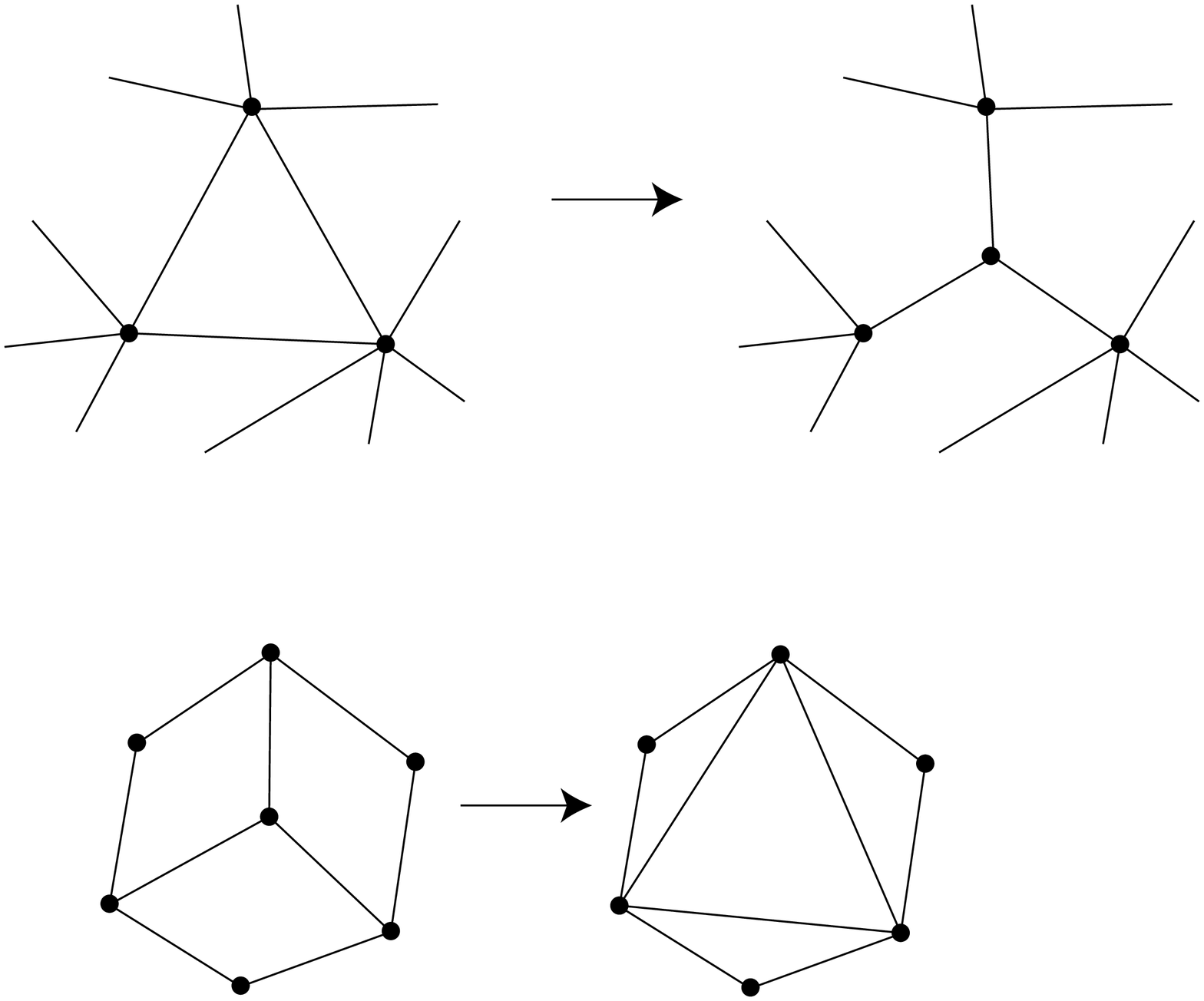}
\end{center}
\vskip-8mm
\caption{\small A $\Delta$-to-$\!Y$ operation.}
\end{figure}

\begin{lemma}\label{lemma:mapG}
  If there is a continuous map
  $f:G\to G'$ that maps disjoint paths and cycles to disjoint paths
  resp.\ cycles,
  and if $G$ is $q$-winding, then $G'$ is also $q$-winding.
\end{lemma}

\begin{proof}
  Let $g:G'\to\R^2$ be any drawing of $G'$. Then $g\circ f:G\to\R^2$
  is a drawing of $G$. Since $G$ is $q$-winding, there are $q$
  disjoint paths or cycles in $G$ that form a $q$-winding partition
  for $g\circ f$. These paths/cycles are mapped under $f$ to $q$
  disjoint paths/cycles in $G'$, which form a $q$-winding partition
  for $g$.
\end{proof}

Any inclusion $G\subset G'$ satisfies the assumptions of this lemma,
as does any series reduction (inverse subdivision of an edge).

\begin{lemma}
  A graph $G'$ obtained from a $\Delta$-to-$\!Y$ operation on a $q$-winding
  graph $G'$ is again $q$-winding. 
\end{lemma}

\begin{proof}
Assume that $G'$ is obtained from $G$ by a $\Delta$-to-$\!Y$ operation, more 
precisely by deleting the edges $v_1v_2,v_2v_3$ and $v_1v_3$ and 
adding the vertex $v$ together with the edges $vv_1,vv_2$ and 
$vv_3$. Define $f:G\to G'$ as the identity on all vertices of $G$ 
and all edges of $G$ except the three deleted ones. For these, 
define $f(v_iv_j):=v_ivv_j$.
The function $f$ maps disjoint paths and cycles to disjoint paths/cycles. 
\end{proof}

\begin{figure}[htb]
\begin{center}
\includegraphics*[bb= 4 260 585 470, scale=.3]{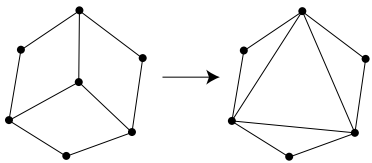}
\end{center}
\vskip-8mm
\caption{\small 
  A $Y\!$-to-$\Delta$ operation that transforms this 2-winding graph into a
  graph that is not 2-winding.} \label{YDeltaCounterexample}
\end{figure}

We note that $Y\!$-to-$\Delta$ operations do not preserve the
property to be $q$-winding: See Figure~\ref{YDeltaCounterexample}.

\begin{lemma}
  If $G$ has a $q$-winding minor, then $G$ is $q$-winding.  
\end{lemma}

\begin{proof}
 The ``minor of'' relation is generated by addition of edges,
 and splitting of vertices. Both operations satisfy the condition
 of Lemma~\ref{lemma:mapG}.  
\end{proof}

We return to the discussion of 2-winding graphs.

\begin{proof}[Proof of Proposition \ref{K_4ANDK_2,3}]
The Winding Number Conjecture holds for  $q=2$, so $K_4$ is 2-winding. The 
graph $K_{2,3}$ is obtained from $K_4$ by a $\Delta$-to-$\!Y$ operation and 
hence is 2-winding  as well.
\end{proof}

\begin{theorem}
A graph is $2$-winding  if and only if it contains $K_4$ or $K_{2,3}$ 
as a minor.
\end{theorem}

\begin{proof}
Every graph that has a $q$-winding  minor is itself $q$-winding. 
Therefore every graph containing $K_4$ or $K_{2,3}$ as a minor is 
2-winding.

On the other hand, if a graph does not contain one of these two 
minors, then it is \emph{outerplanar}, that is, it has 
a planar drawing with all vertices lying on the exterior region.
In such a drawing no two 
edges intersect, and no cycle winds 
around a vertex. Hence the graph is not $2$-winding.
\end{proof}

\subsection{3-Winding  graphs and $q$-winding 
subgraphs of complete graphs}

We prove two general results about $q$-winding subgraphs of
$K_{3q-2}$, and obtain the minimal $3$-winding subgraph of $K_7$.
For the Topological Combinatorics
notation and basics employed in the following,
we refer to Matou\v{s}ek \cite[Chap.~6]{MatousekBZ:BU}.

\begin{theorem}\label{wg1}
Let $p\ge 3$ be a prime and $M$ a maximal matching in 
$K_{3p-2}$. Then $K_{3p-2}-M$ is $p$-winding.
\end{theorem}

\begin{proof}[Proof (suggested by 
Vu{\v{c}}i{\'c} and {\v{Z}}ivaljevi{\'c} {\cite{vuz93},
in the presentation of Matou\v{s}ek \cite[Sect.~6.6]{MatousekBZ:BU}})]~\\
Let $N:=4(p-1)$ and let $f:K_{3p-2}\to\R^2$ be a drawing of 
$K_{3p-2}$, which we may assume to be piecewise
linear in general position.

 We divide the proof in three steps.
\begin{compactenum}[1.~]
\item We describe a $\mathbb{Z}_p$-invariant subcomplex $L$ of 
$(\Delta_N)^{*p}_{\Delta(2)}$. 

\item We show that $\ind_{\mathbb{Z}_p}(L)\ge 
N>N-1=\ind_{\mathbb{Z}_p}((\R^3)^{*p}_{\Delta})$. 
Thus by the defining property of the index
(see Matou\v{s}ek \cite[Sects.~6.2, 6.3]{MatousekBZ:BU}), 
$L$ cannot be mapped to 
$(\R^3)^{*p}_{\Delta}$ $\mathbb{Z}_p$-equivariantly.

\item We extend the drawing $f$ to a map $F:\Delta_N\to\R^3$ and
  examine the Tverberg partitions of~$F$ and winding partitions of $f$
  that are obtained from the equivariant map
  $F^*|_L:L\rightarrow (\R^3)^{*p}$.
\end{compactenum}
\medskip

\noindent
\emph{Step 1:} 
The vertex set of the deleted join complex
$(\Delta_N)^{*p}_{\Delta(2)}$ can be arranged in an array
of $(N+1)\times p$ points, as in Figure~\ref{dots}.
The maximal simplices then have exactly one vertex in each 
of the $N+1$ levels.

We extend the matching $M$ of $K_{3p-2}$ to a maximal 
matching on the vertices of $\Delta_N$ and group the rows into pairs 
accordingly. One row remains single. For each pair of rows we  
choose a $\mathbb{Z}_p$-invariant cycle in the complete 
bipartite graph generated by these two shores, such that the 
cycles contain no vertical edges. (This requires $p\ge3$.)
The maximal simplices of $L$ 
shall be the maximal simplices of $(\Delta_N)^{*p}_{\Delta(2)}$ 
which contain an edge from each cycle (compare Figure~\ref{dots}). 
This defines $L$ as a $\mathbb{Z}_p$-invariant $N$-dimensional
subcomplex of~$(\Delta_N)^{*p}_{\Delta(2)}$.
\medskip

\begin{figure}[ht]
\begin{center}
\includegraphics[bb = 50 550 240 790, scale=.62, clip]{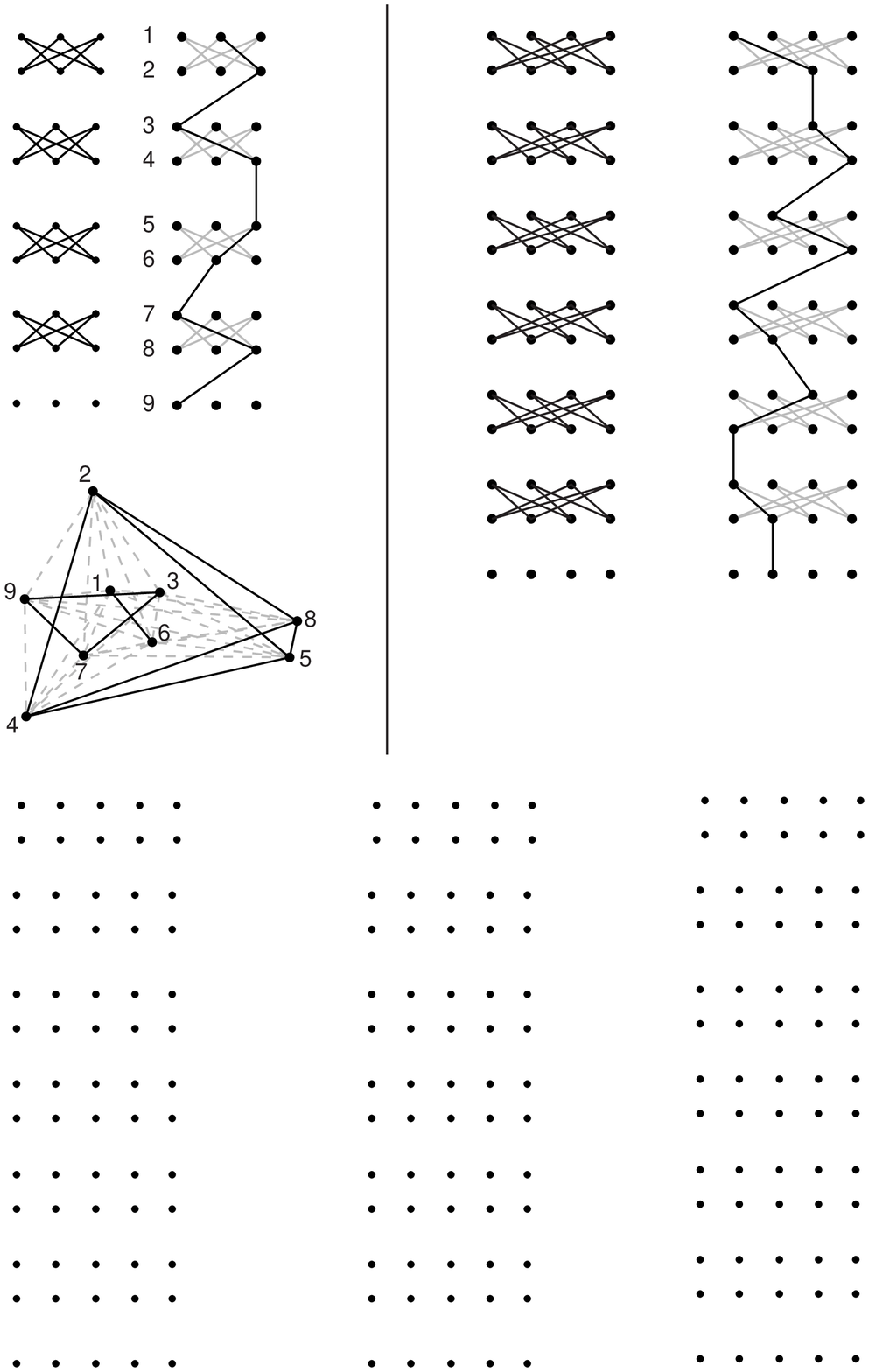}\qquad
\includegraphics[bb = 50 380 240 550, scale=.62, clip]{Dots.eps}
\end{center}
\vskip-5mm
\caption{\small The left figure indicates the 
complex $L$ in the case $p=3$ and $N=8$:
For each pair of rows a cycle is drawn; the 
bold chain indicates a maximal face of $L$.
The figure on the right illustrates the partition of 
the vertices of $(\Delta_8)^{*3}_{\Delta(2)}$ represented by this face.}
\label{dots}
\end{figure}

\noindent
\emph{Step 2:} $L$ can be interpreted as the join of its $N/2$ 
cycles and the remaining row of $p$ points,
\[
L\ \cong\ (S^1)^{*N/2}*D_p\ \cong\ S^{N-1}*D_p.
\]
This space is $N$-dimensional and $(N-1)$-connected, so 
$\ind_{\mathbb{Z}_p}(L)=N$.

The identity $\ind_{\mathbb{Z}_p}((\R^3)^{*p}_{\Delta})=N-1$ 
is elementary as well; see \cite[Sect.~6.3]{MatousekBZ:BU}.
\medskip

\noindent
\emph{Step 3:} Now we can extend $f$ to a continuous map 
$F:\Delta_{4(p-1)}\to\R^3$, such that for every Tverberg partition 
$\{\sigma_1,\dots,\sigma_p\}$ for~$F$, the set 
$\{\sigma_1\cap\Delta_{3(p-1)},\dots,\sigma_p\cap\Delta_{3(p-1)}\}$ 
is a winding partition for $f$. 

According to the pattern of \cite[Sect.~6.3]{MatousekBZ:BU}, this
yields a $\Z_p$-equivariant map 
$F^{*p}:\Delta_{4(p-1)}^{*p}\rightarrow(\R^3)^{*p}$. In view of
the index computation of Step~2, the
restriction $F^{*p}|_L$ hits the diagonal,
which yields a $p$-fold coincidence point in~$L$, and thus
a Tverberg partition for $F:\Delta_{4(p-1)}\to\R^3$
which does not use a matching edge. 
\end{proof}

\begin{proposition}\label{wg2}
Let $X$ be $q-1$ edges of $K_{3q-2}$ meeting in one vertex. Then 
$K_{3q-2}-X$ is not $q$-winding.
\end{proposition}

\begin{figure}
\begin{center}
\includegraphics*[bb= 110 550 470 670,scale=.55]{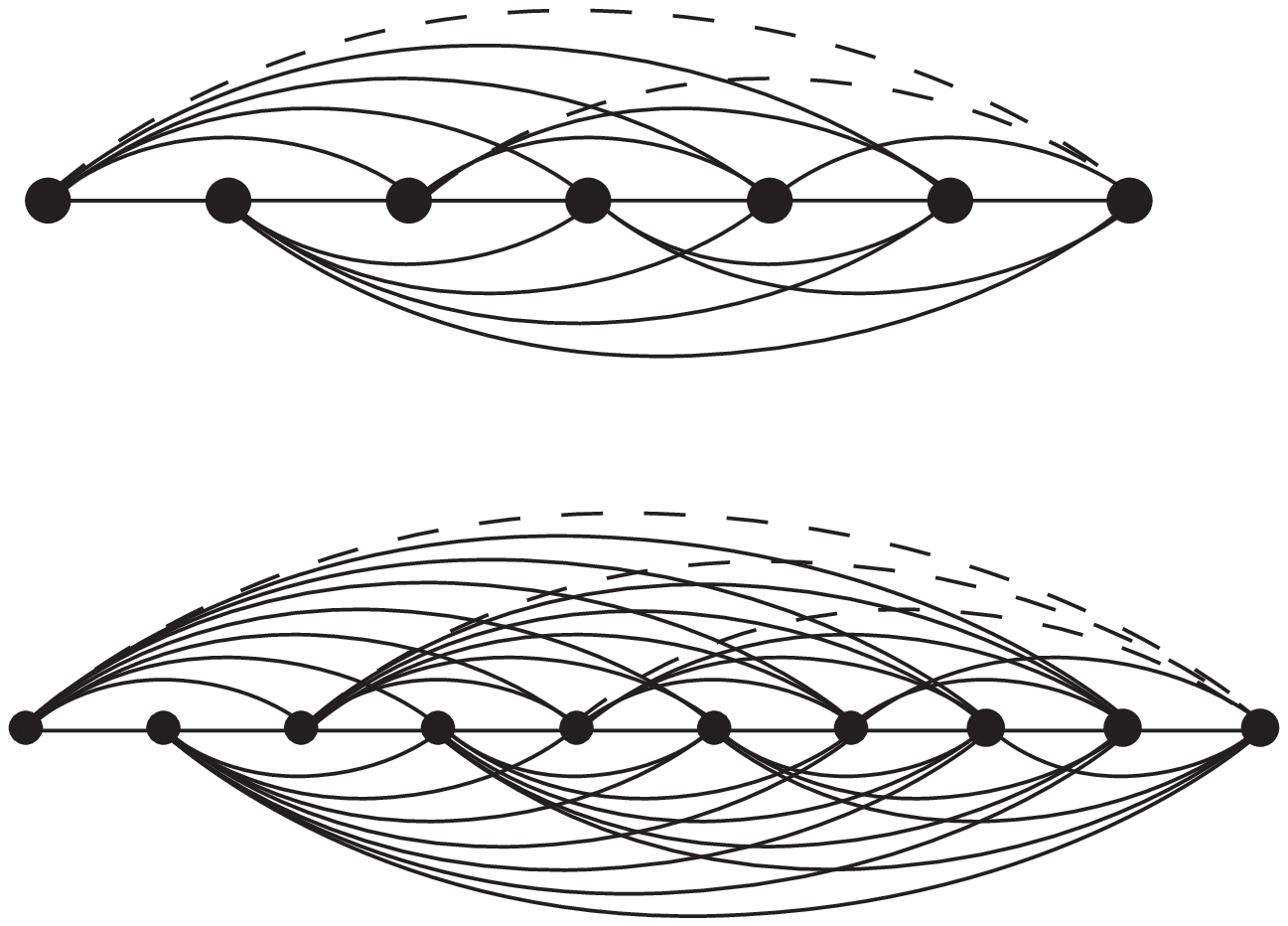}
\includegraphics*[bb= 100 390 500 530,scale=.55]{K_nAlsLinieOhneZweiKanten.eps}
\end{center}
\vskip-3mm
\caption{\small Drawing of $K_7$ and $K_{10}$. The edges that form $X$ 
are dashed.} \label{K_nOhneKante}
\end{figure}

\begin{proof}
All we need is a drawing of $K_{3q-2}-X$ without a 
 $q$-winding partition. We can use the alternating linear model of 
$K_n$ described in Example \ref{K_nByGuy}. 
Order the vertices such that the meeting vertex is at the 
right end of the drawing and the other vertices of $X$ have the 
numbers $1,3,5,\dots,2q-5,2q-3$. The edges of $X$ are then in the 
upper half. (Compare Figure \ref{K_nOhneKante}.)
It is a nice, elementary exercise to verify that in this situation there
is no winding partition that doesn't use an edge of~$X$.
\end{proof}

\begin{corollary}
The unique minimal $3$-winding  minor of $K_7$ is $K_7-M$, where $M$ 
is a maximal matching.
\end{corollary}

\begin{proof}
$K_7-M$ is a 3-winding  minor of $K_7$ (Theorem \ref{wg1}, for $p=3$). It is 
minimal, because all edges not in $M$ are adjacent to an edge in 
$M$ and thus must not be deleted (Proposition \ref{wg2}).

If on the other hand $K$ is a $3$-winding  minor of $K_7$, then only a 
matching can be deleted (again by Proposition \ref{wg2}). For $K$ 
to be minimal, this matching must be maximal.
\end{proof}

\begin{proposition}
Not every $3$-winding  graph has $K_7$ minus a maximal matching as a 
minor.
\end{proposition}

\begin{proof}
Let $M$ be a maximal matching in $K_7$. Execute a $\Delta$-to-$\!Y$ operation 
on~$K_7-M$; the resulting graph is 3-winding, but does not have 
$K_7-M$ as a minor.
\end{proof}
\medskip

\begin{small}
\subsubsection*{Acknowledgements}
Thanks to Stephan Hell and to John Sullivan for interesting
discussions and valuable observations, and to the referees for
several useful comments.

\end{small}
\end{document}

%% file: winding2.pstex_t
\begin{picture}(0,0)%
\includegraphics{winding2.pstex}%
\end{picture}%
\setlength{\unitlength}{1973sp}%
\begingroup\makeatletter\ifx\SetFigFont\undefined%
\gdef\SetFigFont#1#2#3#4#5{%
  \reset@font\fontsize{#1}{#2pt}%
  \fontfamily{#3}\fontseries{#4}\fontshape{#5}%
  \selectfont}%
\fi\endgroup%
\begin{picture}(3460,3462)(2306,-5290)
\put(3976,-5236){\makebox(0,0)[lb]{\smash{{\SetFigFont{9}{10.8}{\familydefault}{\mddefault}{\updefault}{\color[rgb]{0,0,0}$f$}%
}}}}
\put(3526,-2311){\makebox(0,0)[lb]{\smash{{\SetFigFont{9}{10.8}{\familydefault}{\mddefault}{\updefault}{\color[rgb]{0,0,0}$a$}%
}}}}
\put(4426,-2311){\makebox(0,0)[lb]{\smash{{\SetFigFont{9}{10.8}{\familydefault}{\mddefault}{\updefault}{\color[rgb]{0,0,0}$b$}%
}}}}
\end{picture}%

%% file: d=3proof.pstex_t
\begin{picture}(0,0)%
\includegraphics{d=3proof.pstex}%
\end{picture}%
\setlength{\unitlength}{2960sp}%
\begingroup\makeatletter\ifx\SetFigFont\undefined%
\gdef\SetFigFont#1#2#3#4#5{%
  \reset@font\fontsize{#1}{#2pt}%
  \fontfamily{#3}\fontseries{#4}\fontshape{#5}%
  \selectfont}%
\fi\endgroup%
\begin{picture}(4554,3624)(1459,-5173)
\put(1756,-4376){\makebox(0,0)[rb]{\smash{{\SetFigFont{11}{13.2}{\familydefault}{\mddefault}{\updefault}{\color[rgb]{0,0,0}$\wa(f|_{\partial\sigma})$}%
}}}}
\put(2001,-3801){\makebox(0,0)[rb]{\smash{{\SetFigFont{11}{13.2}{\familydefault}{\mddefault}{\updefault}{\color[rgb]{0,0,0}$\partial B_\sigma$}%
}}}}
\put(1936,-2976){\makebox(0,0)[rb]{\smash{{\SetFigFont{11}{13.2}{\familydefault}{\mddefault}{\updefault}{\color[rgb]{0,0,0}$x_i$}%
}}}}
\end{picture}%

%% file: tverberg.bbl
\begin{thebibliography}{10}
\itemsep=0pt

\bibitem{bss81}
{\sc I.~B{\'a}r{\'a}ny, S.~B. Shlosman, and A.~Sz{\H{u}}cs}, {\em On a
  topological generalization of a theorem of {T}verberg}, J. London Math. Soc.
  (2), 23 (1981), pp.~158--164.

\bibitem{deL01}
{\sc M.~de~Longueville}, {\em {N}otes on the topological {T}verberg theorem},
  Discrete Math., 247 (2002), pp.~271--297.
  (The original publication, in a volume of selected papers in honor of
  Helge Tverberg, Discrete Math., 241 (2001), 207--233, 
  was marred by publisher's typography errors.)

\bibitem{hell:_tverb}
{\sc S.~Hell}, {\em On the number of {T}verberg partitions in the prime power
  case}.
\newblock Preprint, TU Berlin, April 2004, 7~pages;
  \url{arXiv:math.CO/0404406}.

\bibitem{MatousekBZ:BU}
{\sc J.~Matou\v{s}ek}, {\em {Using the {Borsuk--Ulam} Theorem. {L}ectures on
  Topological Methods in Combinatorics and Geometry}}, Universitext,
  Springer-Verlag, Heidelberg, 2003.
\newblock Written in cooperation with Anders Bj\"orner and G\"unter M. Ziegler.

\bibitem{Munkres:AT}
{\sc J.~R. Munkres}, {\em Elements of Algebraic Topology}, Addison-Wesley,
  Menlo Park, CA, 1984.

\bibitem{Ozaydin}
{\sc M.~{\"O}zaydin}, {\em Equivariant maps for the symmetric group}.
\newblock Preprint 1987, 17 pages.

\bibitem{saa64}
{\sc T.~L. Saaty}, {\em The minimum number of intersections in complete
  graphs}, Proc. Nat. Acad. Sci. U.S.A., 52 (1964), pp.~688--690.

\bibitem{Sarkaria-primepower}
{\sc K.~S. Sarkaria}, {\em Tverberg partitions and {B}orsuk-{U}lam theorems},
  Pacific J. Math., 196 (2000), pp.~231--241.

\bibitem{dipl-Schoeneborn}
{\sc T.~Sch\"oneborn}, {\em On the Topological {T}verberg Theorem},
  Diplom\-arbeit, 39~pages, TU Berlin 2004; \url{arXiv:math.CO/0405393}.

\bibitem{tve66}
{\sc H.~Tverberg}, {\em A generalization of {R}adon's theorem}, J. London Math.
  Soc., 41 (1966), pp.~123--128.

\bibitem{tve81}
\leavevmode\vrule height 2pt depth -1.6pt width 23pt, {\em A generalization of
  {R}adon's theorem. {II}}, Bull. Austral. Math. Soc., 24 (1981), pp.~321--325.

\bibitem{vol96}
{\sc A.~Y. Volovikov}, {\em On a topological generalization of {T}verberg's
  theorem}, Mat. Zametki, 59 (1996), pp.~454--456.
\newblock Translation in Math. Notes 59 (1996), no. 3-4, 324--325.

\bibitem{vuz93}
{\sc A.~Vu{\v{c}}i{\'c} and R.~T. {\v{Z}}ivaljevi{\'c}}, {\em Note on a
  conjecture of {S}ierksma}, Discrete Comput. Geom., 9 (1993), pp.~339--349.

\bibitem{Z35}
{\sc G.~M. Ziegler}, {\em Lectures on {P}olytopes}, vol.~152 of Graduate Texts
  in Mathematics, Springer-Verlag, New York, 1995.
 Revised edition, 1998.

\end{thebibliography}
